\documentclass[12pt,a4paper]{article}
\usepackage{amssymb}
\title{Isomorphisms of Affine Pl\"ucker Spaces}
\author{Hans Havlicek}
\date{}

\newtheorem{theo}{Theorem}

\newtheorem{remark}{Remark}

\newcommand{\Bcal}{{\cal B}}

\newcommand{\Ecal}{{\cal E}}

\newcommand{\Lcal}{{\cal L}}

\newcommand{\Ncal}{{\cal N}}

\newcommand{\Pcal}{{\cal P}}

\newcommand{\rel}{\sim}
\newcommand{\rrel}{\approx}

\newcommand{\AG}[2]{\mbox{$\mbox{{\rm AG}}(#1,#2)$}}
\newcommand{\AFF}{{\sf A}}
\newcommand{\PL}{\mbox{$(\Pcal,\Lcal,\parallel)$}}
\newcommand{\PLL}{\mbox{$(\Pcal',\Lcal',\parallel')$}}%

\newcommand{\im}{\mbox{\rm im\,}}
\newcommand\inv{^{-1}}

\newcommand{\abb}[3]{\mbox{$#1\,:\,#2\rightarrow#3$}}
\newcommand{\Abb}[5]{\mbox{$#1\,:\,#2\rightarrow#3,\;#4\mapsto #5$}}
\newcommand{\proof}{{\em Proof. }}
\newcommand{\proofend}{$\Box$}
\newcommand{\zitat}[4]{\bibitem{#1}{\sc #2}: {\sl #3\/}. #4.\vspace{-0.0em}}

\begin{document}
\maketitle
   \begin{abstract}
   All isomorphisms of Pl\"ucker spaces on affine spaces with dimensions
   $\geq 3$ arise from collineations of the underlying affine spaces.
   \end{abstract}
\thispagestyle{empty}
\section{Introduction}\label{Intro}

Let $L$ be a set and $\rel$ a reflexive and symmetric binary relation on $L$
such that $(L,\rel)$ is connected, i.e., for any $a,b\in L$ there exists a
finite sequence $a=a_1\rel a_2\rel\cdots\rel a_n=b$. Following {\sc W.\ Benz}
the pair $(L,\rel)$ is called a {\em Pl\"ucker space}\/ \cite[p.\ 199]{Be92}.
Elements $a,b\in L$ are said to be {\em related}\/ if $a\rel b$. {\em
Adjacent}\/ elements ($a\rrel b$) are characterized by $a\rel b$ and $a\neq
b$.

The relation $\not\rrel$ is reflexive and symmetric. However, $(L,\not\rrel)$
is not necessarily a Pl\"ucker space, since it need not be connected.
Nevertheless, $L$ splits into a family of {\em connected components}\/ with
respect to $\not\rrel$, say $(L_i)_{i\in I}$. Each component $L_j $ $(j\in
I)$ gives rise to the Pl\"ucker space $(L_j,{\not\rrel}_j)$, where
${\not\rrel}_j$ denotes the restriction of $\not\rrel$ to $L_j\times L_j$. On
the other hand, $L_j$ is not necessarily connected with respect to $\rel_j$,
i.e., the restriction of $\rel$ to $L_j\times L_j$. Hence $(L_j,\rel_j)$ need
not be a Pl\"ucker space.

Given two Pl\"ucker spaces $(L,\rel)$ and $(L',\rel')$ a bijection
$\abb\varphi {L}{L'}$ is called an {\em isomorphism}\/ if
      \begin{equation}\label{PLUECKER}
      a\rel b \Longleftrightarrow
      a^\varphi\rel' b^\varphi \mbox{ for all } a,b\in L.
      \end{equation}
Obviously, (\ref{PLUECKER}) and
   \begin{equation}\label{PLUECKER_B}
   a\not\rrel b \Longleftrightarrow
   a^\varphi\not\rrel' b^\varphi \mbox{ for all } a,b\in L
   \end{equation}
are equivalent conditions.

All automorphisms of $(L,\rel)$ form its so-called {\em Pl\"ucker group}.
Write, as above, $(L_i)_{i\in I}$ for the connected components of $L$ with
respect to $\not\rrel$. Then each automorphism of $(L_j,{\not\rrel}_j)$
$(j\in I)$ extends to an automorphism of $(L,\rel)$ by setting $x\mapsto x$
for all $x\in L\setminus L_j$.

Let $\AFF=\PL$ be an affine space, where $\Pcal$, $\Lcal$ and $\parallel$
denotes the set of points, the set of lines and the parallelism,
respectively. Lines $a,b\in \Lcal$ are called {\em related}\/ ($a\rel b$), if
$a\cap b\neq\emptyset$. The pair $(\Lcal,\rel)$ is satisfying the conditions
mentioned before and will be called an {\em affine Pl\"ucker space}. We
remark that for $\dim\AFF\geq 3$ the set $\Lcal$ is the set of `points' of
partial linear space, the {\em affine Grassmann space}\/ on $\Lcal$; cf.
\cite{BM82}, \cite{BM82a}, \cite{BM83}, \cite{Ta86} and \cite{Za88}. However,
the relation $\rel$ is not the same as the binary relation of `collinearity'
used in those papers, since `collinear points' are represented by lines that
are related or parallel.

If $\dim\AFF\neq2$, then $(\Lcal,\not\rrel)$ is a Pl\"ucker space. If $\AFF$
is an affine plane, then $(\Lcal,\not\rrel)$ is not a Pl\"ucker space. The
connected components $(\Lcal_i)_{i\in I}$ with respect to $\not\rrel$ are the
pencils of parallel lines, since the relations $\not\rrel$ and $\parallel$
are coinciding now. If $\Lcal_j$ $(j\in I)$ is a fixed pencil of parallel
lines, then the relation ${\not\rrel}_j$ is the coarsest relation on
$\Lcal_j$. Thus Pl\"ucker spaces on affine planes have indeed a very poor
structure. The case $\dim\AFF\leq 1$ cannot deserve interest at all.

We shall determine all isomorphisms of Pl\"ucker spaces on affine spaces
\AFF, $\AFF'$ with dimensions $\geq 3$: Any collineation yields an
isomorphism of the associated affine Pl\"ucker spaces and vice versa. If we
impose additional assumptions on \AFF, $\AFF'$ (cf. Theorem \ref{HINREICH}),
then any bijection $\abb{\varphi}{\Lcal}{\Lcal'}$ is already an isomorphism
of Pl\"ucker spaces whenever (\ref{PLUECKER}) is satisfied with an
implication ($\Longrightarrow$) rather than an equivalence
($\Longleftrightarrow$).

Similar theorems for Pl\"ucker spaces on projective spaces are due to {\sc
W.L.\ Chow} \cite{Ch49}, {\sc H.\ Brau\-ner} \cite{Br88} and the author
\cite{Ha95}. For further results and references on Pl\"ucker spaces see,
among others, \cite{Be92}, \cite{Be94} and \cite{Ha95a}.

\section{Isomorphisms}\label{ISOM}

Let $\AFF=\PL$ and $\AFF'=\PLL$ be affine spaces. If
$\abb\kappa{\Pcal}{\Pcal'}$ is a collineation, i.e., a bijection preserving
collinearity and non-collinearity of points, then $\kappa$ gives rise to a
bijection
   \begin{equation}
   \Abb\varphi{\Lcal}{\Lcal'}{Q\vee R}{Q^\kappa\vee R^\kappa}
   \; (Q,R\in\Pcal,\; Q\neq R)
   \end{equation}
taking related lines to related lines in both directions.

We shall prove the following converse:
   \begin{theo}\label{KAPPA}
   Let\/ $\AFF=\PL$ and\/ $\AFF'=\PLL$ be affine spaces with\/ $\dim\AFF'
   \geq 3$. Suppose that $\abb\varphi\Lcal{\Lcal'}$ is an isomorphism of the
   Pl\"ucker space\/ $(\Lcal,\rel)$ onto the Pl\"ucker space\/
   $(\Lcal',\rel')$. Then
      \begin{equation}\label{DEFKAPPA}
         \Abb\kappa\Pcal{\Pcal'}{a\cap b}{a^\varphi\cap b^\varphi}\;
         (a,b\in\Lcal,\; a\rrel b)
      \end{equation}
      is a well-defined collineation.
   \end{theo}
\proof
{\em (a)}\/
We infer from $\dim\AFF'\geq 3$ and the bijectivity of $\varphi$ that
$\#\Lcal > 1$. Therefore $\dim\AFF\geq 2$. With $Q\in\Pcal$ write $\Lcal(Q)$
for the star of lines with centre $Q$, i.e., the set of all lines in $\Lcal$
running through $Q$. Any star of lines is a maximal set of mutually related
lines.

Suppose that $(\Lcal(Q))^\varphi$ contains a trilateral spanning a plane
$\Ecal'\subset\Pcal'$, say. All lines of $(\Lcal(Q))^\varphi$ are mutually
related. Therefore they are all contained in $\Ecal'$. By $\dim\AFF'\geq 3$,
there exists a line $a\in\Lcal$ with $a^\varphi\cap\Ecal'=\emptyset$. Thus
   \begin{equation}\label{KAPPA1}
   a^\varphi\not\rel' x^\varphi \mbox{ for all }x\in\Lcal(Q)
   \end{equation}
and therefore
   \begin{equation}\label{KAPPA2}
   a\not\rel x \mbox{ for all }x\in\Lcal(Q).
   \end{equation}
On the other hand, there exists a line joining $Q$ with an arbitrarily
chosen point of the line $a$. This contradicts (\ref{KAPPA2}).

Thus we have established that $(\Lcal(Q))^\varphi$ is a subset of a star of
lines for any $Q\in\Pcal$. It is obvious now that (\ref{DEFKAPPA})
is a well-defined mapping.

{\em (b)}\/
Given a point $Q'\in\Pcal'$ one may show as above that
$(\Lcal'(Q'))^{\varphi\inv}$ is a subset of a star of lines. Therefore
$\kappa$ is a surjection and under $\varphi$ stars of lines go over to stars
of lines in both directions.

If points $Q,R\in\Pcal$ are distinct, then
   \begin{equation}
   \#(\Lcal(Q)\cap\Lcal(R))=\#((\Lcal(Q))^\varphi\cap(\Lcal(R))^\varphi)=1
   \end{equation}
and $Q^\kappa\neq R^\kappa$, whence $\kappa$ is injective.

Three mutually distinct points $Q,R,S\in\Pcal$ are collinear if, and only if,
   \begin{equation}
   \#(\Lcal(Q)\cap\Lcal(R)\cap\Lcal(S))=1.
   \end{equation}
This in turn is equivalent to the collinearity of
$Q^\kappa,R^\kappa,S^\kappa\in\Pcal'$. Hence $\kappa$ is a collineation.
\proofend

   \begin{remark}
   {\em
   If the order of $\AFF'$ is greater than two or if $\dim\AFF'\leq2$, then
   any collineation $\Pcal\rightarrow\Pcal'$ is even an {\em affinity}, i.e.
   a collineation preserving parallelism in both directions. Otherwise, the
   existence of a collineation $\Pcal\rightarrow\Pcal'$ implies the existence
   of an affinity $\Pcal\rightarrow\Pcal'$; see \cite[32.5 and 40.4]{Tam69}.
   Hence for $\dim\AFF'\geq 3$ we obtain all isomorphisms of $(\Lcal,\rel)$
   onto $(\Lcal',\rel')$ via the Pl\"ucker group of $(\Lcal,\rel)$ and a
   single affinity $\Pcal\rightarrow\Pcal'$.
   }
   \end{remark}

   \begin{remark}
   {\em
   If $\AFF=\AFF'$, then Theorem \ref{KAPPA} describes the Pl\"ucker group
   for affine spaces with dimension $\geq 3$. This generalizes a result in
   \cite[p.\ 205]{Be92} for real affine spaces%
      \footnote{The proof given there fails to work in case of characteristic
      two.}%
   .
   }
   \end{remark}

   \begin{remark}
   {\em
   Suppose that $\dim\AFF\geq 2$. The following construction yields all
   maximal sets of mutually related lines that are different from stars:
   Choose a point $Q$, an incident line $a$ and a plane $\Ecal$ containing
   $a$. Write $\Lcal(Q,\Ecal)$ for the pencil of lines in $\Ecal$ running
   through $Q$. Next define a family $(\tau_x)_{x\in\Lcal(Q,\Ecal)}$ of
   translations $\abb{\tau_x}{\Ecal}{\Ecal}$ such that
   $\bigcap_{x\in\Lcal(Q,\Ecal)} x^{\tau_x} = \emptyset$. Then
   $\{x^{\tau_x}\mid x\in\Lcal(Q,\Ecal)\}$ is a maximal set of mutually
   related lines other than a star.
   }
   \end{remark}

   \begin{remark}
   {\em
   If $\abb\varphi\Lcal{\Lcal'}$ is an isomorphism and if $\dim\AFF'=2$, then
   $\dim\AFF=2$ according to Theorem \ref{KAPPA}. By virtue of
   (\ref{PLUECKER_B}) it is easy to establish the following result: Pl\"ucker
   spaces on affine planes $\AFF$ and $\AFF'$ are isomorphic if, and only if,
   $\AFF$ and $\AFF'$ have equipotent pencils of parallel lines or, in other
   words, if the order of $\AFF$ equals the order of $\AFF'$.

   The transposition of two distinct parallel lines of an affine plane is an
   example of a Pl\"ucker transformation that does not stem from a
   collineation.
   }
   \end{remark}

We are now going to weaken the assumptions on $\varphi$ in Theorem
   \ref{KAPPA}.

   \begin{theo}\label{UNIDIR}
   Let\/ $\AFF=\PL$ and\/ $\AFF'=\PLL$ be affine spaces with\/ $\dim\AFF'
   \geq 3$. Suppose that $\abb\varphi\Lcal{\Lcal'}$ is a bijection
   satisfying
   \begin{equation}\label{UNIDIR1}
   a\rel b \Longrightarrow
      a^\varphi\rel' b^\varphi \mbox{ for all } a,b\in \Lcal.
      \end{equation}
   Then
      \begin{equation}\label{UNIDIR2}
         \Abb\lambda\Pcal{\Pcal'}{a\cap b}{a^\varphi\cap b^\varphi}\;
         (a,b\in\Lcal,\; a\rrel b)
      \end{equation}
   is a well-defined injection that preserves collinearity and
   non-collinearity of points. Moreover,
   \begin{equation}\label{UNIDIR3}
   \Lcal(Q)^\varphi=\Lcal'(Q^\lambda) \mbox{ for all } Q\in\Pcal.
   \end{equation}
   \end{theo}

\proof
{\em (a)}\/
By the proof of Theorem \ref{KAPPA}, part (a), the following assertions have
been already verified: The dimension of $\AFF$ is $\geq 2$. For all
$Q\in\Pcal$ the set $(\Lcal(Q))^\varphi$ is a subset of a star of lines,
whence (\ref{UNIDIR2}) is a well-defined mapping.

{\em (b)}\/
Let $Q,R\in\Pcal$ be distinct and assume that $Q^\lambda=R^\lambda$. Choose a
line $c\in\Lcal\setminus(\Lcal(Q)\cup\Lcal(R))$ and a point $S\in c$ such
that $Q,R,S$ are not collinear. Therefore $Q\vee S$, $R\vee S$ and $c$ are
three distinct concurrent lines. We deduce from (\ref{UNIDIR2}) that
   \begin{equation}
   S^\lambda = (Q\vee S)^\varphi\cap(R\vee S)^\varphi
   = Q^\lambda = R^\lambda.
   \end{equation}
Consequently, $Q^\lambda\in x^\varphi$ for all $x\in \Lcal$ . This is
impossible due to the surjectivity of $\varphi$. Hence $\lambda$ is
injective.

If $\{Q,R,S\}\subset\Pcal$ is a triangle, then $Q\vee R$, $R\vee S$ and
$S\vee Q$ are three distinct lines. The injectivity of $\varphi$ and the
injectivity of $\lambda$ force that
$\{Q^\lambda,R^\lambda,S^\lambda\}\subset\Pcal'$ is a triangle. By definition,
$\lambda$ is a collinearity-preserving mapping.

{\em (c)}\/
Finally, we establish (\ref{UNIDIR3}). Assume to the contrary that there
exists a point $Q\in\Pcal$ and a line $b\in\Lcal\setminus\Lcal(Q)$ with
$Q^\lambda\in b^\varphi$. Choose two distinct points $R_1,R_2\in b$. Then
$\{Q,R_1,R_2\}$ is a triangle, but
$\{Q^\lambda,R_1^\lambda,R_2^\lambda\}\subset b^\varphi$ is a collinear set,
an absurdity.
\proofend

\vspace{1em}

The aim of the following discussion is to give sufficient conditions for
$\lambda$ to be a collineation or, equivalently, a surjection. We could apply
results on injective mappings of affine spaces preserving collinearity of
points; see \cite[3.1--3.3]{Be92}, \cite{Zi81} and the references in
\cite{Ha94}. However, we proceed instead in close analogy with \cite{Ha95}.
Any star of lines in an affine space, for example, a star $\Lcal(Q)$ in \AFF,
carries in a natural way the structure of a projective space, viz.\
   \begin{equation}
   \AFF/{Q}:=
   (\Lcal(Q),\{\Lcal(Q,\Ecal)\mid Q\in\Ecal, \Ecal\mbox{ a plane}\}).
   \end{equation}
Recall the following concept due to {\sc P.V. Ceccherini} \cite{Ce67}: A {\em
semicollineation}\/ of projective spaces is a bijection taking any three
collinear points to collinear points. The existence of {\em proper}\/
semicollineations (different from collineations) of Desarguesian projective
spaces seems to be an open problem; see also \cite{BT84}, \cite{Kr96} and
\cite{Ma70}. Semicollineations fall within the wider class of {\em weak
linear mappings}\/ that have been characterized independently in \cite{FF94}
and \cite{Ha94}.

Now (\ref{UNIDIR3}) can be improved as follows:
   \begin{theo}\label{SEMI}
   Let\/ $\AFF$, $\AFF'$ and $\varphi$ be given as in Theorem \ref{UNIDIR}.
   Then\/ $\dim\AFF\geq 3$. If moreover the order of\/ $\AFF$ is not two and
   if $Q\in\Pcal$, then the restricted mapping
      \begin{equation}\label{SEMI1}
      \abb{\varphi|\Lcal(Q)}{\Lcal(Q)}{\Lcal'(Q^\lambda)}
      \end{equation}
   is a semicollineation of $\AFF/Q$ onto $\AFF'/Q^\lambda$.
   \end{theo}
\proof
$\AFF$ cannot be the affine plane of order two, since $6<\#\Lcal'$.

Suppose that the order of $\AFF$ is not two. By (\ref{UNIDIR3}), the mapping
(\ref{SEMI1}) is bijective. Let $a,b,c\in\Lcal(Q)$ be `collinear points' of
$\AFF/{Q}$. There exists a line $d\in\Lcal\setminus\Lcal(Q)$ that is adjacent
to $a$, $b$ and $c$. Hence
$a^\varphi,b^\varphi,c^\varphi\in\Lcal(Q^\lambda,Q^\lambda\vee d^\varphi)$
represent `collinear points' in $\AFF'/{Q^\lambda}$ so that (\ref{SEMI1})
is a semicollineation. There are non-coplanar lines through $Q^\lambda$
representing `non-collinear' points of $\AFF'/{Q^\lambda}$. Their pre-images
under (\ref{SEMI1}) are distinct and non-coplanar, so that $\dim\AFF\geq 3$.
\proofend

\vspace{1em}

If $\AFF$ is of order two, then (\ref{SEMI1}) is in general merely a
bijection.

   \begin{theo}\label{HINREICH}
With the settings of Theorem \ref{UNIDIR}, each of the following conditions
is sufficient for $\lambda$ to be a collineation:
\begin{enumerate}
\item $\AFF$ or\/ $\AFF'$ is a finite affine space.
\item $\dim\AFF\leq\dim\AFF'<\infty$.
\item The order of\/ $\AFF$ is different from two and every monomorphism of
      an underlying field $F$ of\/ $\AFF$ in an underlying field $F'$ of\/
      $\AFF'$ is surjective.
\item $\AFF$ and\/ $\AFF'$ are affine spaces of order two.

\end{enumerate}
   \end{theo}
\proof
{\em Ad 1.}
Since $\varphi$ is bijective, both $\AFF$ and $\AFF'$ are finite affine
spaces.

Let $\AFF\cong\AG n2$ and $\AFF'\cong\AG {m}{p^h}$, where $n\geq3$, $m\geq3$,
$h\geq 1$ are integers and $p$ is a prime.
Choose $Q\in\Pcal$. We deduce from (\ref{UNIDIR3}) that
      \begin{equation}
      2^{n-1}\#\Lcal(Q) = \#\Lcal = \#\Lcal' = p^{h(m-1)}\#\Lcal(Q^\lambda) =
      p^{h(m-1)}\#\Lcal(Q).
      \end{equation}
Consequently, $p=2$, $n-1 = h(m-1)$ and, by $\#\Lcal(Q) =
\#\Lcal(Q^\lambda)$,
\begin{equation}\label{SUMME}
   \sum_{i=0}^{n-1} 2^i = \sum_{i=0}^{m-1} 2^{hi}.
   \end{equation}
We infer that each summand on the right hand side of (\ref{SUMME}) appears
exactly once on the left hand side, whence $h=1$ and $n=m$. This implies that
$\lambda$ is surjective.

If the order of $\AFF$ is greater than two, then (\ref{SEMI1}) is a
semicollineation and, by \cite[14.2]{Ce67}, even a collineation. Therefore
\AFF\ and $\AFF'$ are of equal order. Now $\lambda$ turns out to be
   surjective, because of
   \begin{equation}
   \dim\AFF=\dim(\AFF/Q)+1=\dim(\AFF'/Q^\lambda)+1=\dim\AFF'<\infty.
   \end{equation}

{\em Ad 2.}
By virtue of the previous result, we may exclude affine spaces $\AFF$,
$\AFF'$ of order two from the following discussion.

Choose any point $R'\in\Pcal'$. As $\varphi$ is surjective, there exists a
line $d\in\Lcal$ with $R'\in d^\varphi$. Let $Q\in\Pcal$ be off the
line $d$ and put $\Ecal:=Q\vee d$. We observe that
\begin{equation}
   \dim(\AFF/Q) = \dim\AFF-1
   \leq
   \dim\AFF'-1 = \dim(\AFF'/Q^\lambda) < \infty.
   \end{equation}
By \cite[8.4]{Ce67} or \cite[Theorem 2.2]{Kr96} the semicollineation
(\ref{SEMI1}) turns out to be a collineation. Consequently,
$\Lcal(Q,\Ecal)^\varphi$ is a pencil of lines $\Lcal'(Q^\lambda,\Ecal')$,
say. The line $d^\varphi\not\ni Q^\lambda$ is related to all lines of this
pencil with at most one exception. This implies that the point-set
$d^\lambda$ is equal to the affine line $d^\varphi$. Thus
   \begin{equation}\label{DD}
   R'\in d^\varphi= d^\lambda\subset\im\lambda.
   \end{equation}

{\em Ad 3.}
Choose any point $Q\in\Pcal$. By Theorem \ref{SEMI}, the mapping
(\ref{SEMI1}) is a semicollineation. This implies the existence of a
monomorphism $F\rightarrow F'$; cf.\ \cite[5.1]{Ce67}, \cite[Theorem
5.4.1]{FF94} or \cite[Theorem 2]{Ha94}. By \cite[5.3]{Ce67}, the mapping
(\ref{SEMI1}) is a collineation. From this the surjectivity of $\lambda$ is
established as above.

{\em Ad 4.}
Choose any point $R'\in\Pcal'$. Since $\varphi$ is surjective, there exists a
line $\{R_1,R_2\}\in\Lcal$ with
$R'\in\{R_1,R_2\}^\varphi=\{R_1^\lambda,R_2^\lambda\}$.
\proofend

   \begin{remark}
   {\em
   Let $\AFF$ be an affine space over $\mbox{GF(2)}$ with a countable basis
   and let $\AFF'$ be an $m$-dimensional affine space ($3\leq m\leq\aleph_0$)
   over a countable field $F'$ of arbitrary characteristic. Hence there is
   either no monomorphism or no surjective monomorphism of $\mbox{GF(2)}$ in
   $F'$. As $\#\Pcal=\#\Lcal'=\aleph_0$, we can index all points of $\Pcal$
   as $Q_1,Q_2,\ldots$ and all lines of $\Lcal'$ as $a_1',a_2',\ldots$ such
   that there are no repeated elements.

   Let us define, by recursion, an injective sequence
   $\{1,2,\ldots\}\rightarrow\Pcal'$, $s\mapsto R_s'$ such that each line of
   $\Lcal'$ contains exactly two points: We start with a point $R_1'\in
   a_1'$ and put $\Bcal'_1:=\{R_1'\}$. Next assume that we are already given a
   set $\Bcal_i'=\{R_1',\ldots,R_i'\}$ formed by $i\geq 1$ mutually distinct
   points no three of which are collinear. Write $\Ncal_i'$ for the set of
   all lines that arise by joining distinct points of $\Bcal_i'$. Then let
   $j\in\{1,2,\ldots\}$ be the least element such that the line $a_j'$ is not
   in $\Ncal_i'$. Since $a_j'$ carries an infinite number of points, we can
   choose such a point $R_{i+1}'\in a_j'\setminus\Bcal_i'$ that no three
   elements of the set $\Bcal_{i+1}':=\Bcal_i'\cup\{R_{i+1}'\}$ are
   collinear.

   Put $\Bcal':=\bigcup_{s=1}^\infty \Bcal'_s$. By construction, no three
   distinct points of $\Bcal'$ are collinear. Furthermore, given a line
   $a_k'\in\Lcal'$ we obtain that $a_k'\in\Ncal'_{2k}$, as required.

   The mapping
   \begin{equation}
   \Abb{\varphi}{\Lcal}{\Lcal'}{\{Q_s,Q_t\}}{R'_s\vee R'_t},\quad
   (s,t\in\{1,2,\ldots\},\;s\neq t)
   \end{equation}
   is a bijection satisfying (\ref{UNIDIR1}). The associated
   injection $\lambda$ (see (\ref{UNIDIR2})) takes $Q_s$ to $R'_s$
   ($s\in\{1,2,\ldots\}$). Only two points of $a_1'$ belong to $\im\lambda$,
   whence $\lambda$ is not surjective.
   }
   \end{remark}

   \begin{remark}{\em
   Let $\AFF$, $\AFF'$ be affine spaces with equal infinite order and
   $2=\dim\AFF'<\dim\AFF\leq\aleph_0$. There exists a bijection
   \abb{\varphi}{\Lcal}{\Lcal'} such that any class of parallel lines in
   $\AFF$ is mapped onto a pencil of parallel lines in $\AFF'$. Such a
   $\varphi$ is satisfying (\ref{UNIDIR1}) without being an isomorphism of
   Pl\"ucker spaces.
   }\end{remark}

   \begin{remark}
   {\em
   Let $\AFF=\AFF'$ be an affine plane of infinite order. Choose two
   non-parallel lines $a,b\in\Lcal$. There exists a bijection
   \abb{\varphi}{\Lcal}{\Lcal} such that the parallel class of $a$ is mapped
   onto the union of the parallel classes of $a$ and $b$, whereas any other
   pencil of parallel lines is mapped onto a pencil of parallel lines. Then
   $\varphi$ is satisfying (\ref{UNIDIR1}) without being a Pl\"ucker
   transformation.
   }\end{remark}

\vspace {1em}

\noindent
Hans Havlicek\\
Abteilung f\"ur Lineare Algebra und Geometrie\\
Technische Universit\"at\\
Wiedner Hauptstra{\ss}e 8--10\\
A-1040 Wien, Austria\\
EMAIL: {\tt havlicek@geometrie.tuwien.ac.at}

\end{document}